\documentclass[12pt,leqno]{article}
\usepackage{a4wide,epsf,latexsym}
\usepackage[centertags]{amsmath}
\usepackage{graphicx}

\newcommand{\eps}{\varepsilon}

\begin{document}

\title {Monotone discretization of elliptic problems with mixed derivatives on anisotropic meshes:
a counterexample}

\author{Hans-G. Roos, TU Dresden}

\date{January 2019}

\maketitle

In many papers of O'Riordan and Shishkin the authors underline the importance of numerical approximations
for singularly perturbed elliptic problems which are free of spurious oscillations. Therefore they prefer
(difference) methods where the associated system matrix is a monotone matrix, and they use, exclusively,
the special class of M-matrices, which are monotone. But what about elliptic problems with
mixed derivatives?

Not much is known concerning finite difference methods for singularly perturbed elliptic problems
with {\it mixed derivatives on layer-adapted meshes}, in \cite{RST96} and  \cite{SS09} is also nothing to find. It is
well-known that on isotropic meshes one can generate an M-matrix (see, for instance, Theorem 10.1 in
\cite{Ma02}). That means, that a condition of the type
\[
      C_1\le \frac{h_x}{h_y}\le C_2
\]
is sufficient to generate an M-matrix. But layer-adapted meshes are highly anisotropic. In \cite{HO17}
the authors state that a monotone scheme imposes a condition on the mesh ratio, but this is not
proved so far. In \cite{DOS09} the authors avoid to discuss this question: they use simple an inconsistent approximation
which approximates the mixed derivative only on the fine, isotropic part of a Shishkin mesh. The
assumption $\eps\le N^{-1}$ then allows nevertheless to prove an error estimate.

{\it In this paper we present an example which shows that a consistent first order scheme for an elliptic
problem with mixed derivatives on a highly anisotropic mesh cannot generate an M-matrix.}

Consider the elliptic operator
\[
   Lu:=u_{xx}+u_{xy}+u_{yy}
\]
on the rectangle $(x_0-H,x_0+H)\times(y_0-h,y_0+h)$ and the nine-point difference operator
\begin{subequations}
  \begin{align}
  L_{H,h}u  & :=\alpha u(x_0,y_0)+\beta_1 u(x_0-H,y_0)+\beta_2 u(x_0+H,y_0)+\beta_3 u(x_0,y_0-h)+\beta_4 u(x_0,y_0+h)\nonumber\\                          &+\beta_5 u(x_0-H,y_0-h)+\beta_6 u(x_0+H,y_0-h)+\beta_7 u(x_0-H,y_0+h)+\beta_8 u(x_0+H,y_0+h).\nonumber
   \end{align}
\end{subequations}

Taylor expansion yields the first order consistency conditions
\begin{equation}\label{0}
\alpha+\sum \beta_i=0
\end{equation}
\begin{equation}\label{1}
\beta_2-\beta_1+\beta_6+\beta_8-(\beta_5+\beta_7)=0
\end{equation}
\begin{equation}\label{2}
\beta_4-\beta_3+\beta_7+\beta_8-(\beta_5+\beta_6)=0
\end{equation}
\begin{equation}\label{3}
\frac{H^2}{2}(\beta_1+\beta_2+\beta_5+\beta_6+\beta_7+\beta_8)=1
\end{equation}
\begin{equation}\label{4}
\frac{h^2}{2}(\beta_3+\beta_4+\beta_5+\beta_6+\beta_7+\beta_8)=1
\end{equation}
and
\begin{equation}\label{5}
hH(\beta_5-\beta_7+\beta_8-\beta_6)=1.
\end{equation}
Fixing $\beta_2, \beta_3, \beta_4$, the parameter $\beta_1$ is given by
\[
  H^2(\beta_1+\beta_2)=h^2(\beta_3+\beta_4),
\]
and it is easy to show that the remaining parameters are uniquely determined.
\vspace{0.2cm}

{\it The question is: Can one choose the parameters in such a way that $\beta_i\ge 0$
for all $i$ ?}

\vspace{0.3cm}
\noindent
In the isotropic case this is possible, for instance, for $h=H$ one gets with the choice
$\beta_2=\beta_3=\beta_4=\frac{1}{2h^2}$ the result $\beta_1=\beta_5=\beta_7=\frac{1}{2h^2}$, moreover
$\beta_6=\beta_8=0$ and $\alpha=-\frac{3}{2h^2}$. The corresponding matrix is the negative of
an M-matrix.

Now let us assume that \eqref{1}--\eqref{5} do have a solution with $\beta_i\ge 0$ for all $i$ in
the case $h<<H$.

From \eqref{3} we obtain $0\le \beta_i\le \frac{2}{H^2}$ for $i=1,2,5,6,7,8$. Consequently, \eqref{2}
implies
\begin{equation}
\beta_4-\beta_3={\cal O}(\frac{1}{H^2}).
\end{equation}
The equations \eqref{2} and \eqref{5} yield a system for $B:=\beta_6-\beta_8$ and
$D:=\beta_5-\beta_7$:
\begin{equation}
B+D=\beta_4-\beta_3, \quad -B+D=\frac{1}{hH}.
\end{equation}
Thus we get $B=-\frac{1}{2hH}+\frac{1}{2}(\beta_4-\beta_3)$, respectively
\begin{equation}
\beta_6=-\frac{1}{2hH}+\beta_8+\frac{1}{2}(\beta_4-\beta_3).
\end{equation}
Because $\beta_8$ and $\beta_4-\beta_3$ are of order ${\cal O}(\frac{1}{H^2})$, they cannot compensate
the large negative term $-\frac{1}{2hH}$ if $h$ is sufficiently small in comparison to a given $H$.
Consequently, $\beta_6$ cannot be nonnegative for  anisotropic meshes where $H/h$ is sufficiently
large, and this is the case for Shishkin meshes and other types of layer-adapted meshes.

{\it Author's address}: Hans-Goerg Roos,\\ Faculty of Mathematics, Technical University Dresden,
 Zellescher Weg 12-14,\\ 01062 Dresden, Germany
  \\e-mail: Hans-Goerg.Roos@tu-dresden.de

\end{document}